\documentclass[11pt]{amsart}
\RequirePackage[left=1in,right=1in,top=1in,bottom=1in]{geometry}
\usepackage{amssymb}
\usepackage{amsfonts}
\usepackage{amscd}
\usepackage{amsmath}
\usepackage{mathrsfs}
\usepackage{curves}
\usepackage{epsfig}
\usepackage{graphicx}
\usepackage{verbatim}
\usepackage{latexsym}
\usepackage{eucal}
\usepackage{hyperref}
\usepackage{comment}

\usepackage{tikz}
\usepackage{tikz-3dplot}

\hypersetup{
    colorlinks,
    citecolor=blue,
    filecolor=magenta,
    linkcolor=red,
    urlcolor=black
}

\numberwithin{equation}{section}

\newtheorem{theorem}{Theorem}[section]

\newtheorem{definition}[theorem]{Definition}
\newtheorem{remark}[theorem]{Remark}
\theoremstyle{definition}

\def \Rm {\mathbb R}

\def \Sm {\mathbb S}

\newcommand{\dint}{\displaystyle\int}

\def \mC {{\mathcal C}}

\def \beqq {\begin{equation}}
\def \eeqq {\end{equation}}

\def \bpf {\begin{proof}}
\def \epf {\end{proof}}
\def \beq {\begin{equation*}}
\def \eeq {\end{equation*}}

\def \eps {\epsilon}  

\def \tilde {\widetilde}

\def \red {\color{red}}
\def \fu {\mathfrak{u}}

\newcommand{\gb}[1]{{{\color{blue}{GB: #1}}}}
\newcommand{\bp}[1]{{{\color{magenta}{BP: #1}}}}

\begin{document}
\title{Fermi pencil beams and Off-axis Laser Detection}
\author{Guillaume Bal$^\dagger$}
\thanks{$^\dagger$Departments of Statistics and Mathematics, University of Chicago, Chicago, IL.\\ {\em Email address:} \texttt{guillaumebal@uchicago.edu}}
\author{Benjamin Palacios$^\ddagger$}
\thanks{$^\ddagger$Department of Mathematics, Pontificia Universidad Cat\'olica de Chile, Santiago, Chile.
\\ {\em Email address:} \texttt{benjamin.palacios@mat.uc.cl}}
\maketitle


%
%
%
%
%
%
%
%
%
%
%
%
\begin{abstract}
    This paper concerns the reconstruction of properties of narrow laser beams propagating in turbulent atmospheres. We consider the setting of off-axis measurements, based on light detection away from the main path of the beam. 
    
    We first model light propagation in the beam itself by macroscopic approximations of radiative transfer equations that take the form of Fermi pencil beam or fractional Fermi pencil beam equations. Such models are effective in the small mean-free-path, large transport-mean-free path regime. The reconstruction of their constitutive parameters is also greatly simplified compared to the more accurate radiative transfer equations or (fractional) Fokker-Planck models.

    From off-axis measurements based on wide-angle single scattering off the beam, we propose a framework allowing us to reconstruct the main features of the beam, and in particular its direction of propagation {\em and} the location of the emitting source.
\end{abstract}

\section{Introduction}
This paper concerns the reconstruction of a narrow laser beam propagating in a turbulent atmosphere; see \cite{BG18,C-SPIE-03,HBDH,RR-OE-08} for references and applications. We assume that the turbulence has correlation length and the laser a central wavelength that are both very small compared to the overall distance of propagation. In such a setting, light propagation is accurately modeled by a deterministic radiative transfer (transport) equation \cite{BKR,BKR-liouv}, as used in \cite{C-SPIE-03}; see \eqref{eq:tr}-\eqref{eq:Q} below. 

Although inverse transport theory is well developed (see, e.g., \cite{B} and references there), forward models based on the kinetic transport involve a high dimensional reconstruction of its constitutive parameters. In the presence of limited available data, it is advisable to look for more macroscopic models to describe beam propagation. 

In the highly forward-peaked regime, when light interacts often with the underlying turbulence (small mean free path) but at each scattering event with small variation in its direction (larger transport mean free path), two types of equations are known to emerge: Fokker-Planck and fractional Fokker-Planck models; see \eqref{eq:FP} and \eqref{eq:fFP} below. The derivation of such Fokker-Planck models from radiative transfer equations in the highly forward-peaked regime is done, e.g., in \cite{AS,BL1,BL2,P}. Inverse transport theory (i.e., the reconstruction of their constitutive coefficients from boundary measurements) is mathematically open.

It turns out that in the regime of large transport mean-free path (small diffusion coefficient) another accurate approximation is possible. It is based on neglecting back-scattering and is thus valid as long as beams remain sufficiently narrow. The derivation of the Fermi pencil beam and fractional Fermi pencil beam models from Fokker-Planck models was derived recently in \cite{BP-SIMA-20,BP-preprint-20}. We present the models in detail in section \ref{sec:fpb}.

Inverse problems based on the Fermi pencil beam models are significantly simpler than those based on radiative transfer or Fokker-Planck equations. Moreover, they offer a parametric (the `fraction' coefficient) set of models for beam spreading based on the possibly unknown statistics of turbulence and thus provide reasonable macroscopic descriptions for the reconstruction of laser beams from off-axis measurements.

\medskip

Section \ref{sec:kmodels} recalls the aforementioned kinetic models (radiative transfer, Fokker-Planck) and presents the setting for off-axis measurements: a decomposition of light scattering into two components, $Q_F$ modeling forward-peaked beam scattering, and $Q_S$ modeling small, large-angle scattering that generates the signals captured by off-axis detectors. 

As mentioned above, the Fermi pencil beam models are presented in section \ref{sec:fpb} in detail. We also show their accuracy in a suitable metric when compared to (fractional) Fokker-Planck solutions in the small diffusion regime.

The settings for the off-axis measurements are presented in section \ref{sec:rec}. We provide a sufficient set of measurements such that the main axis of propagation of the laser may be estimated. Assuming a cylindrical symmetry of the laser beam, we propose an inversion based on an inverse Radon transform that explicitly provides reconstruction for the turbulence diffusion (including the `fraction' coefficient) as well as the location of its source from a minimal set of off-axis measurements.

\section{Kinetic model}\label{sec:kmodels}
We model light propagation in scattering media with a kinetic equation of the form
\begin{equation}\label{eq:tr}
  \theta\cdot\nabla_x w + \lambda w = Q(w) + f,
\end{equation}
where $x\in X\subset \Rm^n$ is spatial position (with $n=3$ in most applications), $\theta\in\Sm^{n-1}$ the angular direction and $w(x,\theta)$ the particle density. We assume a light speed normalized to $1$ and a constant index of refraction. The parameter $\lambda(x)\geq0$ models intrinsic absorption while $Q(w)$ is a scattering kernel. In the radiative transfer model, the scattering kernel is of the form
\begin{equation}\label{eq:Q}
  Q(w) = \dint_{\Sm^{n-1}} k(x,\theta,\theta') (w(x,\theta')-w(x,\theta)) d\theta',
\end{equation}
with $d\theta$ the standard Lebesgue measure on the sphere and $k(x,\theta,\theta')$ quantifying scattering from $\theta'$ to $\theta$ at position $x$. 

Our objective is to model the propagation of narrow beams of light and to address their reconstruction from off-axis measurements, i.e., measurements performed away from the physical location of the main beam. Beams preserve their structure in environments with large transport mean free path, heuristically defined as a distance over which light direction changes significantly. Thus, for a distance between $\theta$ and $\theta'$ sufficiently large, we expect the scattering kernel $k$ to be small. In many regimes of interest, the mean free path, defined as the (average) distance between successive interactions of light with the underlying medium, may still be quite small. This is the regime of forward-peaked scattering, where $k$ is significant for $\theta$ close to $\theta'$. We thus arbitrarily separate $k$ into a component with $\theta$ close to $\theta'$ and a component $k_S$ where this is not the case.  We correspondingly decompose 
\[
 Q = Q_F + Q_S
\]
with $Q_F$ corresponding to forward-peaked scattering while $Q_S$ models small, large-angle scattering. 

Since $Q_S$ is small, we perform an expansion in that parameter and define
\[
   \theta\cdot\nabla_x u + \lambda u= Q_F(u) + f,\qquad  \theta\cdot\nabla_x w + \lambda w = Q_F(w) + Q_S(u) + f.
\]
The first equation for $u$ is the main object of interest to analyze beam spreading. The second equation for $w$ is an approximation of the original $w$ to second order in $Q_S$. Off-axis scattering and measurements are then modeled by $u_S=w-u$. Neglecting spreading $Q_F(u_S)$ and absorption $\lambda u_S$ in the scattered contribution to simplify reconstructions, which is justified for instance when the distance between the beam and the detector arrays is reasonably short, then $u_S$ is the solution of
\begin{equation}\label{offaxis_dist}
   \theta\cdot\nabla_x u_S
   = Q_S(u).
\end{equation}

We now focus on the different models for $Q_F(u)$. The scattering kernel \eqref{eq:Q} may be derived from models of wave propagation in heterogeneous media when the wavelength of the wave packets is comparable to the correlation length of the random medium \cite{BKR}. When the correlation length is significantly larger than the wavelength, the limiting equation is of the form of a Fokker-Planck equation instead \cite{P,BKR}, with 
\[
  Q_F= Q_{FP} :=  \nabla_\theta\cdot D(x,\theta) \nabla_\theta 
\]
where $D(x,\theta)$ is a positive tensor. To simplify, we assume that $D$ is isotropic (independent of $\theta)$ and $Q_F$ then takes the standard form of a Laplace (Beltrami) operator on the unit sphere:
\begin{equation} \label{eq:FP}
Q_F= D(x) \Delta_\theta.
\end{equation}
While the Fokker-Planck equation may be formally derived from the radiative transfer model, it was shown in \cite{P} and derived rigorously in \cite{AS,GPR} that fractional Fokker-Planck models may be better macroscopic approximations for highly forward-peaked regimes. More precisely and following \cite{AS}, consider scattering kernels of the (generalized) Henyey-Greenstein form
\beq\label{HG_kernels}
k_g(x,\theta,\theta') =  \frac{D(x)}{\big(\frac12(1-g)^2+g(1-\theta'\cdot\theta)\big)^{\frac{n-1}{2}+s}}.
\eeq
The scattering kernel is isotropic (depends only on $\theta\cdot\theta'$) and is forward-peaked when $g<1$ is close to $1$. The standard Henyey-Greenstein phase function is obtained for $n=3$ and $s=\frac12$. We consider a parameter $0<s<1$. The  limiting scattering kernel as $g\to1$ formally takes the form
\[
   k_F(x,\theta,\theta') =   \frac{D(x)}{\big(1-\theta'\cdot\theta\big)^{\frac{n-1}{2}+s}}.
\]
We then verify that the mean free path vanishes since $\int_{\Sm^{n-1}} k_F(x,\theta,\theta')d\theta'=+\infty$. However, $k_F(x,\theta,\theta')\big(u(\theta')-u(\theta)\big)$ {\em is} integrable for $u$ sufficiently regular and the transport mean free path is finite. Beam structures then appear when the latter is actually large. We will call a fractional Fokker Planck equation the limiting kinetic model with scattering operator
\begin{equation}\label{eq:fFP}
  Q_F u (\theta)=Q_{fFP} u(\theta) :=  \dint_{\Sm^{n-1}} \frac{D(x)}{\big(1-\theta'\cdot\theta\big)^{\frac{n-1}{2}+s}} \big(u(\theta')-u(\theta)\big)d\theta'.
\end{equation}
We refer to \cite{AS,GPR} for a rigorous convergence of radiative transfer solutions to the above fractional Fokker-Planck model. 

We thus obtain a family of fractional Fokker-Planck equations for $0<s<1$ and a standard Fokker-Planck equation formally corresponding to the case $s=1$. It remains to consider the regime of propagation where beam structures may be observed. We obviously need a source term $f$ concentrated in phase space in the vicinity of a point $(x_0,\theta_0)$. We then need to ensure that scattering does not significantly alter the initial directional information. This imposes that the transport mean free path (the mean free path vanishes in all Fokker-Planck models) be large compared to the distance of propagation of interest. 

\begin{figure}
\centerline{
\includegraphics[scale=0.2]{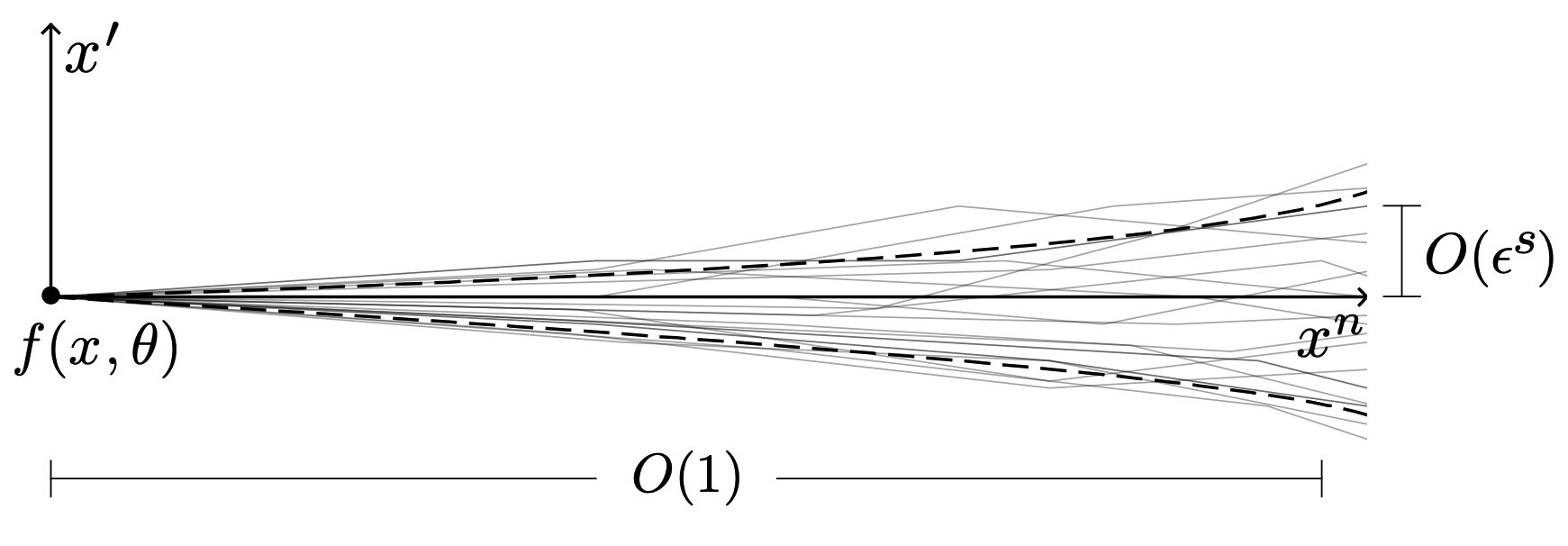}
}
\caption{Spreading of the Fokker--Planck solution}
\label{fig:beam}
\end{figure}

Consider again a transport model of the form
\[
  \theta\cdot\nabla_X u - Qu + \lambda u=0.
\]
The short distance problem consists of defining $X=\eps^{2s} x$ (for some $s\in(0,1]$ depending on the physics of the problem) and then recasting the above as
\[
  \theta\cdot\nabla_x u + \eps^{2s} Q u + \eps^{2s} \lambda_\eps u=0.
\] 
We then assume that $\lambda=\eps^{2s} \lambda_\eps>0$ is a small but leading order term. This ensures that dissipation is large enough to prevent larger-scale phenomena to perturb the analysis. This regime of large transport mean free path is the right one to preserve beam structures and leads (following the decomposition presented above) to the small diffusion (fractional) Fokker-Planck problem (see Fig. \ref{fig:beam}).

We observe that scattering is a perturbation of ballistic transport formally obtained when $\eps=0$. The latter is a reasonable approximation of the beam location but fails to account for any dispersion (beam spreading). It turns out that a higher-order expansion in $\eps$ yields the (fractional) Fermi pencil beam models and that these models accurately capture beam dispersion when $\eps$ is small.

Summarizing the above derivations, the solutions $u$ describing the propagation of a beam generated by the source term $f$ and $u_S$ describing off-axis measurements satisfy the following coupled system:
\begin{equation}\label{FP_and_offaxis}
   \left\{\begin{array}{l}
   \theta\cdot\nabla_x u + \eps^{2s} Q_{F} u + \lambda u=f,\\[2mm]
   \quad\theta\cdot\nabla_x u_S= Q_S(u).
   \end{array}
   \right.
\end{equation}

%
%
%
%
%
%

%
\section{Fermi pencil beam approximation} 
\label{sec:fpb}
%


Neglecting nonlinear effects (which may be important) here, (fractional) Fokker-Planck models are appropriate models to describe laser beam propagation in turbulent atmospheres as we recalled in the preceding section. 
Solving Fokker-Planck equations remains however challenging, both theoretically and computationally. As mentioned above, the simplest approximation is the free (ballistic) transport solution of
\[
\theta\cdot \nabla_x u + \lambda u = 0.
\]
Ballistic models have been used used previously to tackle the off-axis laser detection problem; see, for instance, \cite{HBDH,HB}. However, they ignore important information related to the broadening of the beam which is crucial, for instance, to determine parameters such as the source location. We thus propose to use a more accurate approximation of the Fokker-Planck equation which takes the form of a Fermi pencil beam model.



Our starting point for the beam model is a Fokker-Planck (FP) or fractional Fokker-Planck (fFP) equation with small diffusion given by
\begin{equation}\label{u_fFP}
\theta\cdot\nabla_x u + \lambda u + \epsilon^{2s}D(x)\mathcal{I}^s_\theta[u] = f,\quad (x,\theta)\in\Rm^n\times\Sm^{n-1},
\end{equation}
for a narrow source term $f$ concentrated in the vicinity of a phase space point $(0,\vec{e}_n)\in\Rm^n\times\Sm^{n-1}$. To cover the local (FP) and non-local (fFP) cases simultaneously, we introduce the notation $\mathcal{I}^s_\theta[u]$ to represent the Laplace-Beltrami operator on the unit sphere $-\Delta_\theta $ when $s=1$, while for $s\in(0,1)$ we set  $\mathcal{I}^s_\theta = (-\tilde{\Delta}_\theta)^su - cu$ with constant $c=\frac{\Gamma(\frac{n-1}{2}+s)}{\Gamma(\frac{n-1}{2}-s)}>0$ (here $\Gamma$ stands for the Gamma function), and $(-\tilde{\Delta}_\theta)^s$ defined in stereographic coordinates $v=\mathcal{S}(\theta)$ as the following version of the Laplacian
\begin{equation}\label{fracLapB}
[(-\tilde{\Delta}_\theta)^su]_{\mathcal{S}} := \frac{1}{2^{2s}}\langle v\rangle^{n-1+2s}(-\Delta_v)^s\left(\frac{[u]_\mathcal{S}}{\langle \cdot\rangle^{n-1-2s}}\right).
\end{equation}
We now define the terms that appear in \eqref{fracLapB}. The stereographic coordinates and the associated surface measure are defined as (see \cite[p. 35]{Lee}) $\mathcal{S}:\Sm^{n-1}\backslash\{(0,\dots,0,-1)\}\to \Rm^{n-1}$ where
\[
\begin{aligned}
&v = \mathcal{S}(\theta):= \frac{1}{(1+\theta_n)}(\theta_1,\dots,\theta_{n-1}),
\quad \text{and}\quad d\theta=\frac{2^{n-1}}{\langle v\rangle^{2(n-1)}}dv,\quad\text{with}\quad \langle v\rangle = (1+|v|^2)^{1/2};
\end{aligned}
\]
while the inverse stereographic transformation is defined as
\[
\theta = \mathcal{S}^{-1}(v) := \left(\frac{2v}{\langle v\rangle^2}, \frac{1-|v|^2}{\langle v\rangle^2}\right).
\]
The term $[u]_\mathcal{S}$ corresponds to the particle density $u$ in stereographic coordinates. Moreover, $(-\Delta_v)^s$ stands for the standard (Euclidean) fractional Laplacian given by the singular integral
\[
(-\Delta_v)^sg(v) := c_{n-1,s}\;\text{p.v.}\int_{\Rm^{n-1}}\frac{g(v) - g(v+z)}{|z|^{n-1+2s}}dz,\quad \text{for}\quad s\in(0,1),
\]
for a constant $c_{n-1,s}^{-1}=\int_{\Rm^n}\frac{1-e^{i\hat{\xi}\cdot z}}{|z|^{n+2s}}dz>0$.  We refer the reader to \cite{AS} for details on this version of the Laplace-Beltrami operator.
\begin{remark}
The above definition of the Laplacian differs by a factor $2^{-2s}$ from the one used by the authors in \cite{BP-preprint-20}. It allows for a consistent normalization of the diffusion coefficient when passing to stretched coordinates and deducing the Fermi pencil-beam equation.
\end{remark}

\begin{figure}
\centerline{
\includegraphics[scale=0.25]{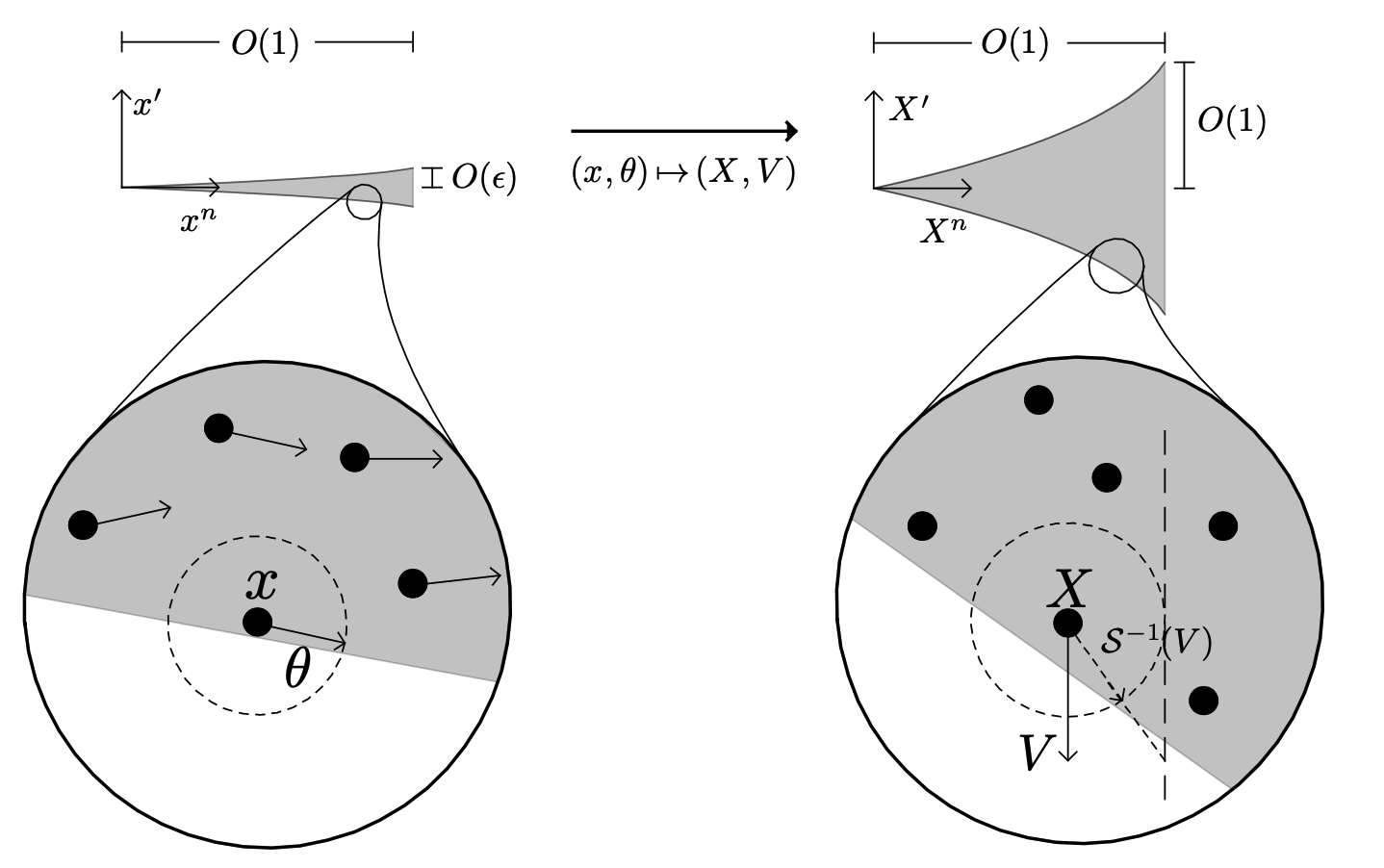}
}
\caption{Narrow beam in pencil-beam (or stretched) coordinates.}
\label{fig:coordinates}
\end{figure}

The above diffusion coefficient is scaled such that diffusion away from the main axis of the beam $\{t(0,\dots,0,1),\ t\geq0\}$, occurs at the scale $\eps$ in phase space. The pencil-beam approximation is based on neglecting backscattering, which is justified in narrow beams with small diffusion $\eps\ll1$. Such a diffusion is naturally captured by the 
%
%
%
following {\em pencil-beam coordinates} (or {\em stretched} coordinates):
\begin{equation}\label{pb_coordinates}
X = ((2\epsilon)^{-1}x',x^n) \quad\text{and}\quad V = \epsilon^{-1}\mathcal{S}(\theta),\quad (x,\theta)\in\Rm^n\times\Sm^{n-1}.
\end{equation}
See Figure \ref{fig:coordinates} for the geometry of the pencil-beam coordinates. 


\begin{definition}\label{def:Fpb}
We say $U(X,V)$ is a (fractional) pencil-beam if it solves the (fractional) Fermi pencil-beam equation (FBP, respectively fFPB) 
\begin{equation}\label{fFPB}
\partial_{X^n}U + V\cdot \nabla_{X'}U + \tilde{\lambda} U + \tilde{D} (-\Delta_V)^sU = 0,\quad (X,V)\in\Rm^n_+\times\Rm^{n-1}
\end{equation}
for $s=1$ (respectively, $s\in(0,1)$), and with a singular boundary source
\begin{equation}\label{fFPB_source}
U(X',0,V) = F_0\delta(X')\delta(V).
\end{equation}
For any $s\in(0,1]$, it takes the explicit form (in the Fourier domain) 
\begin{equation}\label{eq:solFPB}
\mathcal{F}_{X',V}[U](\xi,X^n,\eta) = F_0 e^{-\int^{X^n}_0\tilde{\lambda}(r)dr} e^{-\int^{X^n}_0|\eta+(X^n-t)\xi|^{2s}\tilde{D}(t)dt}.
\end{equation}
\end{definition}

The coefficients $\tilde{\lambda}$ and $\tilde{D}$ are functions of $X^n$ only and determine the attenuation and broadening of the pencil-beam along the main axis. They are related to the coefficients $\lambda$ and $D$ of the Fokker-Planck model as follows: denoting by $\vec{e}_n = (0,\dots,0,1)\in \Rm^n$ and the map $X\mapsto\tilde{X}:=X^n\vec{e}_n$, then  
\begin{equation}\label{eq:tilde}
\tilde{\lambda}(X) := \lambda(\tilde{X})\quad \text{and}\quad \tilde{D}(X):= \frac{1}{2^{2s}}D(\tilde{X}) \quad\text{for}\quad s\in(0,1].
\end{equation}
Note that we use $\;\tilde{\;}\;$ to restrict coordinates or coefficients to the $X^n$-axis (the beam's axis).

A pencil-beam solution may be written in the following self-similar form:
\begin{equation}\label{Uself-sim}
U(X,V) = c_{n-1}\frac{F_0}{(X^n)^{n-1+\frac{n-1}{s}}}\mathfrak{J}\left(\frac{X'}{(X^n)^{1+\frac{1}{2s}}},X^n,\frac{V}{(X^n)^{\frac{1}{2s}}}\right)\text{exp}\left(-\int^{X^n}_0\tilde{\lambda}(t)dt\right)
\end{equation}
for some appropriate constant $c_{n-1}>0$ and with $\mathcal{J}$ defined in the Fourier domain by
\[
\mathcal{F}_{X',V}[\mathfrak{J}](\xi,X^n,\eta):=\text{exp}\left(-\int^1_0|\eta+t\xi|^{2s}\tilde{D}(X^n(1-t))dt\right).
\]
For more details on this and the properties of $\mathfrak{J}$, see \cite{BP-preprint-20}. 

To quantify the accuracy of the ballistic and FPB approximations, we consider a metric that penalizes by {\em how far} particles are in the approximate model compared to where they should be in a FP model. More precisely, we have:
\begin{definition}
Given to two positive Radon measures $f,g$, their {\em $(1,\kappa)$-Wasserstein distance} is given by
$$
\mathcal{W}^1_{\kappa}(f,g):=\sup\left\{\int_{\Rm^n\times \Sm^{n-1}}\psi(f-g) : \psi \text{ Lipschitz with }\|\psi\|_\infty\leq 1,\; \text{Lip}(\psi)\leq \kappa\right\}.
$$
\end{definition}
The gauge of how far particles are from where they should be is given in units of $\kappa^{-1}$ in the sense that $\mathcal{W}^1_{\kappa}(\delta_x,\delta_y)=\kappa|y-x|$. The metric captures beam broadening in the sense that  $\mathcal{W}^1_{\kappa}(\varphi_{\epsilon_1},\varphi_{\epsilon_2})$ is proportional to $\kappa|\epsilon_1-\epsilon_2|$, for $\varphi_{\epsilon_i}(x) = \frac{1}{\epsilon_i^{n}}\varphi\left(\frac{x}{\epsilon_i}\right)$, $i=1,2$, two approximations to the delta function.

We assume that the source term $f\in L^\infty_{x,\theta}\cap L^1_{x.\theta}$ is highly concentrated around $(0,\vec{e}_n)\in\Rm^n\times\Sm^{n-1}$ and more procisely that 
\begin{itemize}
\item[(a)] $f\in L^\infty_{x,\theta}\cap L^1_{x.\theta}$ is compactly supported and for some $\delta>0$ small, 
\[\left| \int f\varphi dxd\theta - F_0\varphi(0,\vec{e}_n)\right|\lesssim \delta  \qquad \mbox{for all $\varphi \in C(\Rm^n\times\Sm^{n-1})$}.
\]
\end{itemize}

In what follows, we denote by $u$ the solution to the FP or fFP equation in \eqref{u_fFP} with a source as above.
We denote by $v$ the solution to the ballistic transport:
\begin{equation}\label{eq:ballistic}
\theta\cdot\nabla_x v + \lambda v = f,\quad (x,\theta)\in\Rm^n\times\Sm^{n-1};
\end{equation}
 while the pencil-beam approximation is defined as
\begin{equation}\label{Fpb}
\mathfrak{u}(x,\theta) = \frac{H(x^n)}{(2\epsilon)^{2(n-1)}}U((2\epsilon)^{-1}x',x^n,\epsilon^{-1}\mathcal{S}(\theta)),
\end{equation}
with $H$ the Heaviside step function and $U$ solution to the FPB equation when $s=1$, or the fFPB equation otherwise. This is the pullback of $U$ with respect the coordinate transformation \eqref{pb_coordinates}, extended by zero to $x^n<0$, and amplified by the factor $(2\epsilon)^{-2(n-1)}$ to preserve the $L^1$-norm to leading order.

The next result (see \cite[Theorem 1.1]{BP-SIMA-20} and \cite[Theorem 1.1]{BP-preprint-20}) summarizes the approximation errors between the various beam models in terms of the small diffusion magnitude $\epsilon>0$ and the resolution parameter $\kappa$.

\begin{theorem}\label{thm:W_est}
Assume that $f$ satisfies (a) above with $\delta\lesssim \kappa \epsilon^{2s}$. 
For a fractional exponent $s\in(0,1)$, we have that for any $s'\in (0,s)$ in dimension $n\geq 3$, or any $s'\in(2s-1,s)$ in dimension $n=2$, there exist positive constants $A(n,s),B(n,s,s')$ and $C(n,s,s')$ such that
\begin{equation}\label{ineq_thm}
A\min\{\kappa\epsilon,1\}\leq\mathcal{W}^1_\kappa(u,v)\leq B (\kappa\epsilon)^{\min\{2s',1\}} \quad\text{and}\quad
 \mathcal{W}^1_\kappa(u,\mathfrak{u})\leq C \kappa^{s'}\epsilon^{2s'},
\end{equation}
where $B\to\infty$ as $s'\to s \leq 1/2$, otherwise, $B$ is independent of $s'$, and $C\to\infty$ as $s'\to s$.

For the local case $s=1$, there are positive constants $A(n),B(n),C(n)$ such that
\[
\mathcal{W}^1_\kappa(u,v)\leq B\kappa\epsilon \quad\text{and}\quad \mathcal{W}^1_\kappa(u,\mathfrak{u})\leq C \kappa\epsilon^2,
\]
and for $\kappa \gtrsim \epsilon^{-1}$ we also have
$A\kappa\epsilon \leq \mathcal{W}^1_\kappa(u,v)$.
\end{theorem}

The above result states that the Fermi-pencil beam model $\mathfrak{u}$ is always a more accurate approximation of the Fokker-Planck solution $u$ than the ballistic model $v$. Moreover, when $\kappa=\eps^{-1}$, i.e., when errors in the location of the particles are gauged in the natural scale $\eps$ of beam spreading, then the ballistic transport is inaccurate, as is obvious physically, while the Fermi-pencil beam model retains a reasonable accuracy (of order $\eps$ when $s=1$ for instance).



From an inversion perspective, the fractional Fermi-pencil beam model has several advantages: it accurately models beam spreading, which ballistic approximations do not (and hence cannot possibly be used to reconstruct the source location) and at the same time has a reasonably explicit expression as recalled in \eqref{eq:solFPB}, which is not the case of the more accurate (fractional) Fokker-Planck model. Leaving the fraction $s$ unknown also provides an additional parameter to model the statistical properties of the (unknown) turbulence. The next section on the off-axis reconstruction problem uses the Fermi-pencil beam model to describe beam spreading.

\section{Beam parameter reconstructions} 
\label{sec:rec}
\subsection{Off-axis measurements}
%
%

Off-axis measurements are modeled following the decomposition $w=u+u_S$ introduced in section \ref{sec:kmodels}, with $u$ and $u_S$ solutions to \eqref{FP_and_offaxis} with $f\geq0$. In this decomposition, $u$ models the beam's particle density which we can approximate, for instance, by letting $u=\fu$ the pencil-beam solution of \eqref{Fpb}, 
while $u_S$ represents the off-axis contribution, thus satisfying
\begin{equation}\label{off_axis_light}
\theta\cdot\nabla_x
u_S(x,\theta) = Q_S(u):=\int_{\Sm^2} \sigma(x,\zeta,\theta)u(x,\zeta)d\zeta.
\end{equation}
For the rest of the section, we assume 
that $\sigma=\sigma(x)$ is isotropic to simplify the analysis.
We use a system of coordinates with origin the source of the laser and $\vec{e}_3=(0,0,1)$ the main direction of the beam. We define $x=(x^1,x^2,x^3)\in\Rm^3$ with $x'=(x^1,x^2)\in \Rm^2$, and similarly $\theta=(\theta^1,\theta^2,\theta^3)\in \Sm^2$.
We then deduce that at $(x,\theta)$ the density of photons takes the explicit form
\begin{equation*}
u_S(x,\theta)=\int^\infty_0 \int_{\Sm^2}
\sigma(x-t\theta)u(x-t\theta,\zeta)d\zeta dt.
\end{equation*}

A flat screen (array) optical camera is modeled as follows. For $x_0$ the center of the array, $L_0$ its fixed radius (with arrays assumed circular for concreteness), and $\theta_0\in \Sm^2$ its direction, we assume that a camera measures light intensity on the set $\mathcal{C}_{x_0,\theta_0}=\{(x,\theta_0):(x-x_0)\cdot\theta_0=0,\; |x-x_0|<L_0\}$, i.e., a disk orthogonal to $\theta_0$ of radius $L_0$ centered at $x_0$.

In this idealized model, an array measures light that comes (exactly) orthogonally to the screen. We assume that the camera can be arbitrarily rotated around the point $x_0$. We also assume available measurements for camera centers $x_0\in X\subset\Rm^3$. We therefore assume measurements known for 
\begin{equation}\label{eq:measdata}
  \Rm^3\times \Sm^2 \supset \Sigma:= \cup_{x_0\in X} \{(x,\theta)\in \Rm^3\times \Sm^2; \ (x-x_0)\cdot\theta=0,\ |x-x_0|< L_0\}.
\end{equation}
The off-axis measurements are therefore characterized by $m(x,\theta)=u_S(x,-\theta)$ for $(x,\theta)\in\Sigma$.

We define the measurement operator $\mathcal{M}$, mapping the unknown parameters to the available information, as
\begin{equation*}
\mathcal{M}:(f,\lambda,D,s,\sigma)
\longmapsto m(x,\theta) := u_S(x,-\theta)|_{\Sigma},
\end{equation*}
for measurements given explicitly by the integrals
\begin{equation}\label{meas1}
m(x,\theta)= \int^\infty_0 \int_{\Sm^2} 
\sigma(x+t\theta)u(x+t\theta,\zeta)d\zeta dt\quad\forall (x,\theta)\in \Sigma.
\end{equation}

Note that the system of coordinates used above, and in particular its origin at the source location, remain unknown.


Our capability to recover any of the parameters $(f,\lambda,D,s,\sigma)$ depends on the available measurement set $\Sigma$ and on any prior knowledge we may possess on the parameters. We now present analytical reconstructions of some parameters of the problem. We start with the determination of the beam's axis and then address the inverse problem of determining the location of the laser's source under suitable assumptions.

\subsection{Determining the beam's axis}
A triangulation procedure is first used to determine an approximation of the central axis of the beam. The beam has a spatial width of order $\eps\ll1$. At this level of approximation, we may model the measurements using the ballistic model $v$ solution of 
$\theta\cdot\nabla v =f_\eps$ with $f_\eps$
a source term of the form 
$f_\eps(x,\theta) = \varphi_{\epsilon}(x)\delta_{\vec{e}_3}(\theta)$, with $\varphi_{\epsilon}(x)=\frac{1}{\epsilon^3}\varphi(\epsilon^{-1}x)$, $\varphi$ a smooth nonnegative and compactly supported function near $x=0$, and where $\delta_{\vec{e}_3}$ stands for the Dirac delta on the unit sphere with support $\{\vec{e}_3\}$. 

Since $v$ is given explicitly by
\[
v(x,\theta)=\int^\infty_0 
f_\eps(x-t\theta,\theta)dt = \int^\infty_0
\varphi_{\epsilon}(x-t\theta)\delta_{\vec{e}_3}(\theta)dt,
\]
the observations (which follows by replacing $u$ with $v$ in \eqref{meas1}) take the form
\[
\begin{aligned}
m(x,\theta) &= \int^\infty_0 \int_{\Sm^2} \int^\infty_0\sigma(x+t\theta)
\varphi_{\epsilon}(x+t\theta-\tau\zeta)\delta_{\vec{e}_3}(\zeta)d\tau d\zeta dt \\
&= \int^\infty_0  \int^\infty_0\sigma(x+t\theta)
\varphi_{\epsilon}(x+t\theta-\tau \vec{e}_3)d\tau dt.
\end{aligned}
\]

\begin{figure}
\begin{center}
\tdplotsetmaincoords{65}{110}
\begin{tikzpicture}[tdplot_main_coords,font=\sffamily]
\fill[red] (3,0,0) circle (1.5pt); \fill[red] (2,0,0) circle (1.5pt); \fill[red] (1,0,0) circle (2pt);
\fill[red] (0,0,0) circle (2.5pt); \fill[red] (-1,0,0) circle (3pt);\fill[red] (-2,0,0) circle (3.7pt);
\fill[red](-3,0,0) circle (5pt);
\draw[thick] (3,0,0) -- (-3,0,0);
\draw[fill=red,opacity=0.3] (3,-3,-2) -- (3,3,2) -- (-3,3,2) -- (-3,-3,-2) -- cycle;
\draw[thick] (3,0,0) -- (-3,0,0);
\draw[fill=blue,opacity=0.3] (-3,-3,0) -- (-3,3,0) -- (3,3,0) -- (3,-3,0) -- cycle;

\draw[fill=green,opacity=0.5] (-.5+1,-3,-.5) --  (-.5+1,-3,.5) -- (.5+1,-3,.5) -- (.5+1,-3,-.5) --cycle;
\draw[thick,red] (.5,-3,0) -- (1.5,-3,0);
\draw[fill=green,opacity=0.5] (-.5-2,-3,-.5) --  (-.5-2,-3,.5) -- (.5-2,-3,.5) -- (.5-2,-3,-.5) --cycle;
\draw[ultra thick,red] (-2.5,-3,0) -- (-1.5,-3,0);
\draw[fill=green,opacity=0.2] (.35,-3.1,-1.8) --  (-.35,-3.1,-1.8)   -- (-.35,-2.9,-2.2) -- (.35,-2.9,-2.2) --cycle;
\draw[very thick,red] (.35,-3,-2) -- (-.35,-3,-2);
\node[anchor=south west,align=center] (line) at (8,7,2.5) {Source};
\draw[-latex] (line) to[out=180,in=-75] (3.05,0.05,-0.05);
\node[anchor=south west,align=center] (line) at (0,4,1.9) {Beam spreading};
\draw[-latex] (line) to[out=180,in=75] (-3.05,0.05,0.15);
\end{tikzpicture}
\end{center}
\caption{Geometry of measurements: Red: beam spreading; green: three detectors used for triangulation and beam spreading along (now approximately known) axis of propagation.}
\label{fig:geom}
\end{figure}
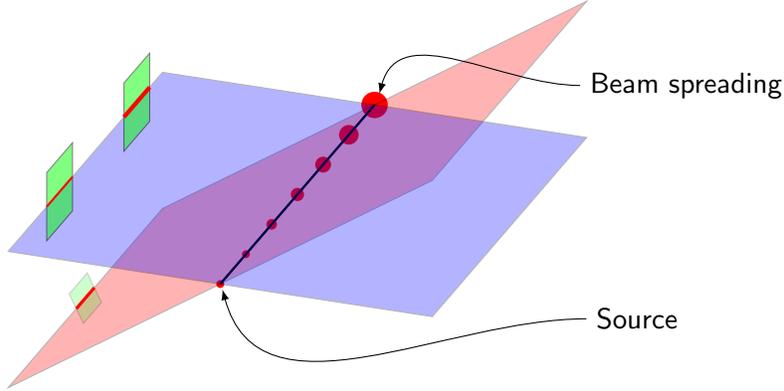

Consider now a camera centered at $x_0$. We orientate the camera so as to maximize the intensity measured at $(x_0,\theta_0)$ for $\theta=\theta_0$. For this choice of $\theta$, we observe that the measurements provide another direction $\phi_0\in \Sm^2$ such that the measurements $m(x_0+t\phi_0,\theta_0)$ are maximal in the sense that $m(x_0+t\phi_0+\delta \theta_0\times \phi_0,\theta_0)$ decay rapidly in $|\delta|$. This requires that the detector array be sufficiently large $L_0\gg\eps$; see Fig. \ref{fig:geom} for an illustration.

With the above information, we obtain that the main direction of the beam belongs to the plane defined by $(x_0,\theta_0,\phi_0)$. 

Let us now assume the existence of other detectors at $x_j$ for $1\leq j\leq J$ and therefore other planes $(x_j,\theta_j,\phi_j)$ constructed as above. Then, the main direction of the beam belongs to the intersection of these planes. We thus require a minimum of two detectors at $x_0$ and $x_1$ such that the corresponding planes are different and hence intersect along the main direction of the beam;  see Fig. \ref{fig:geom}.

The above simple triangulation procedure provides an approximate location of the beam axis with an error proportional to the beam width in the absence of additional information on its structure. As indicated above, it relies on (i) a proper orientation of the camera; and (ii) a sufficiently large array (or aperture) $L_0$ compared to the beam's width to identify the direction $\phi_j$. 

For a determination of the beam's axis based on stochastic properties of the light measurements, which are entirely neglected in the deterministic models considered in this paper, see \cite{BG18}.

\subsection{Beam structure measurements}\label{sec:beam_meas}
We obtained in the preceding section an approximation of the main direction of the beam. The objective of this section is to use the (fractional) Fermi pencil-beam models to reconstruct additional features of the beam including its source location. We therefore assume measurements given by \eqref{meas1} with $u=\fu$ a Fermi pencil-beam solution.

We assume now that $\lambda$, $D$ and $\sigma$ are {\em constant}. In particular, $\tilde{\lambda}=\lambda$ and $\tilde D=\frac{1}{2^{2s}}D$ in \eqref{eq:tilde}; see definition \ref{def:Fpb}.
However, we do not assume them to be known. We also assume here the source term to be a delta function at $x=0$ and direction $\theta=\vec{e}_3$. To relate the measurements $m(x,\theta)$ in \eqref{meas1} for some $(x,\theta)\in \Sigma$, to the unknown coefficients, we first notice that
\[
\begin{aligned}
m(x,\theta) &= \sigma\int^\infty_0 \int_{\Sm^2} 
\fu(x+t\theta,\zeta)d\zeta dt\\
&= \sigma\int^\infty_0 
\int_{\Sm^2}\frac{1}{(2\epsilon)^{4}}U((2\epsilon)^{-1}(x'+t\theta'),x^3+t\theta^3,\epsilon^{-1}\mathcal{S}(\zeta))d\zeta dt\\
&= \sigma\int^\infty_0
\int_{\Rm^2}\frac{1}{(2\epsilon)^{2}}U((2\epsilon)^{-1}(x'+t\theta'),x^3+t\theta^3,V)\langle \epsilon V\rangle^{-2(3-1)}dV dt.
\end{aligned}
\]
For measurements taken at a direction perpendicular to the beam (i.e., $\theta\cdot \vec{e}_3=0$), 
and approximating $\langle \epsilon V\rangle^{-4}$ by $1$, we obtain that
\[
\begin{aligned}
 m(x,\theta) & = \sigma\int^\infty_0
\int_{\Rm^2}\frac{1}{(2\epsilon)^2} U((2\epsilon)^{-1}(x'+t\theta'),x^3,V)dVdt + E(\epsilon), 
\\
E(\epsilon) = E(\epsilon;x,\theta) &:=  \sigma\int^\infty_0
\int_{\Rm^2}\frac{1}{(2\epsilon)^2}U((2\epsilon)^{-1}(x'+t\theta'),x^3,V)\left(\langle \epsilon V\rangle^{-4}-1\right)dV dt.
\end{aligned}
\]

Following \eqref{eq:solFPB} with $n=3$, we compute
\[
\begin{aligned}
m(x,\theta)-E(\epsilon)&=\sigma\int^\infty_0
\int_{\Rm^2}\frac{1}{(2\epsilon)^2} U((2\epsilon)^{-1}(x'+t\theta'),x^3,V)dVdt\\
&=\sigma\int^\infty_0\frac{
1}{(2\epsilon)^{2}} \mathcal{F}^{-1}_{X'}\left[\mathcal{F}_{X',V}[U](\xi,x^3,0)\right] \left((2\epsilon)^{-1}(x'+t\theta')\right)dt\\
&= \sigma F_0\int^\infty_0\frac{e^{-\lambda
x^3}}{(2\epsilon)^{2}}  \mathcal{F}^{-1}_{X'}\left[e^{-2^{-2s}D|\xi|^{2s}\int^{x^3}_0|x^3-z|^{2s}dz}\right] \left((2\epsilon)^{-1}(x'+t\theta')\right) dt\\
&= \sigma F_0\int^\infty_0 e^{-\lambda
x^3} \mathcal{F}^{-1}_{X'}\left[e^{-\epsilon^{2s}D|\xi|^{2s}\int^{x^3}_0|x^3-z|^{2s}dz}\right] \left(x'+t\theta'\right) dt,
\end{aligned}
\]
where the last line follows from the scaling properties of the Fourier Transform: $\mathcal{F}[f(\delta x)](\xi) = \delta^{-2}\mathcal{F}[f(x)](\delta^{-1}\xi)$. 
Therefore,
\begin{equation}\label{meas3}
m(x,\theta)= \int^\infty_0 
g(x+t\theta) dt + E(\epsilon)
\end{equation}
with
\begin{equation}\label{meas32}
g(x):= \sigma F_0e^{-\lambda x^3}\mathcal{F}^{-1}_{X'}\left[e^{-A_{s}(x^3)|\xi|^{2s}}\right] \left(x'\right), \qquad  A_{s}(x^3) = \epsilon^{2s}D\int^{x^3}_0|x^3-z|^{2s} dz.
\end{equation}

We now estimate $E(\epsilon)$.
After a change of variables $t\to t/(2\eps)$ in the above expression, we observe that 
\[
\begin{aligned}
m(x,\theta) &
= \frac{\sigma}{2\epsilon}\int^\infty_0
\int_{\Rm^2}U((2\epsilon)^{-1}x'+t\theta',x^3,V)\langle \epsilon V\rangle^{-4}dV dt.
\end{aligned}
\]
Let us decompose $x=x\cdot\theta\theta+x\cdot\theta^\perp\theta^\perp+x^3\vec{e}_3$ with $\theta^\perp=\vec{e}_3\times\theta$. A measurement array $\mC$ is parametrized by $(x\cdot\theta^\perp,x^3)$ with a support of order $O(\epsilon)$ in the first variable. Since the integral of $m(\cdot,\theta)$ over $\mC$ is independent of $\epsilon$, we observe that $\|m(\cdot,\theta)\|_{L^\infty(\mC)}=O(\epsilon^{-1})$.
The relative error between our measurements and the line integrals $\int^\infty_0 g(x+t\theta) dt$ is thus given by
\[
\frac{\|m(\cdot,\theta) - \int^\infty_0 g(\cdot+t\theta) dt\|_{L^\infty(\mathcal{C})}}{\|m(\cdot,\theta)\|_{L^\infty(\mathcal{C})}} \lesssim \|\epsilon E(\epsilon)\|_{L^\infty(\mathcal{C})}.
\]
From the definition of $E(\epsilon)$ and a Taylor expansion of $\langle \epsilon V\rangle^{-4}-1$, we find that
\[
|2\eps E(\epsilon)|\leq C_{s'} \epsilon^{2s'} \sigma \int^\infty_0
\int U((2\epsilon)^{-1}x'+t\theta',x^3,V)|V|^{2s'}dVdt,
\]
for some  constant $C_{s'}>0$. 
We now recall \eqref{Uself-sim}. From similar estimates on $\mathcal{F}^{-1}[e^{-t|\xi|^{2s}}](x)$ (see \cite{BG}), we deduce the decay properties of the fundamental solution $\mathfrak{J}$: for a fixed $X^n>0$ there are constants $C>c>0$ such that
\[
\frac{c}{(1+|X'|^2+|V|^2)^{\frac{1}{2}(n-1+2s)}}\leq \mathfrak{J}(X',V;X^n)\leq \frac{C}{(1+|X'|^2+|V|^2)^{\frac{1}{2}(n-1+2s)}}.
\]
Then
\[
|2\eps E(\epsilon)|\leq C_{x^3,\sigma,s'}\epsilon^{2s'}\int^\infty_0\int_{\Rm^2}\frac{
|V|^{2s'}}{(1+|(2\epsilon)^{-1}x'+t\theta|^2+|V|^2)^{\frac{1}{2}(n-1+2s)}}dVdt, 
\]
while for the local case of $s=1$ (Gaussian beam), the decay of $\mathfrak{J}$ is exponential. 
Consequently, since $|x'|\gg\epsilon$ (measurements are taken far away from the beam), the integrals above are bounded independently of $\epsilon$ and we conclude that 
\[
\frac{\|m(\cdot,\theta) - \int^\infty_0 g(\cdot+t\theta) dt\|_{L^\infty(\mathcal{C})}}{\|m(\cdot,\theta)\|_{L^\infty(\mathcal{C})}}  = O(\epsilon^{2s'})
\]
for any $0<s'<s<1$ (and in fact $s'=1$ if $s=1$ as one may verify). The constant in the error estimate blows up as $s'$ approaches $s\in(0,1)$.

Notice that $g(x',x^3)=g(|x'|,x^3)$ is radial in the transverse variable.  Since measurements are obtained away from the beam, we observe that $m(x,\theta)$ gives, to leading order in $\epsilon$, the integral $\int_{\ell}
g(y)d\ell(y)$ with $\ell$ the line passing through $x$ in the direction $\theta$. 
The previous analysis then shows that the available measurements provide information on the constitutive coefficients of the (f)FPB models up to an error that is consistent with the approximation of (f)FP by (f)FPB.

\medskip

We now describe a measurement setting allowing us to reconstruct $g(x)$ for suitable values of $x$. 
In the above simplified setting, the beam has a rotational symmetry in $x'$.  Consider a camera centered at $x_0$. We choose the orientation of the camera $\theta_0$ such that the measurements at $x_0+r\phi_0+t\psi_0$ for $\psi_0=\phi_0\times\theta_0$ are symmetric in $t\to-t$ and such that $\theta_0$ is orthogonal to the estimated direction of the beam. Up to an error proportional to $\eps$, the camera is orientated as depicted in Fig.\ref{fig:geom} and $\phi_0=\vec{e}_3$ the main direction of the beam. We therefore have access to the measurement $m(x_0+t\phi_0+\tau \psi_0,\theta_0)$ for $t^2+\tau^2<L_0^2$. Let us fix $t$ with $|t|$ sufficiently small. Then $\tau\mapsto m(x_0+t\phi_0+\tau \psi_0,\theta_0)$ provides (an $\eps-$ approximation of) the line integral
\[
 \tau\mapsto R_1g(\tau)= \int_{\Rm}
 g(x_0+t\phi_0+\tau\psi_0+\mu\theta_0)d\mu.
\]
Since $L_0\gg\epsilon$, we may assume that $R_1g(\tau)=0$ for large values of $\tau$.

For a fixed value of $t$ corresponding to a fixed value of the coordinate $x^3=x_0^3+t\phi_0^3$, we therefore obtain the line integrals (for the two-dimensional set of lines orthogonal to $\phi_0=\vec{e}_3$) of the function $x'\mapsto g(x',x^3)=g(|x'|,x^3)$ since line integrals in other directions $\theta$ such that $\theta\cdot \phi_0=0$ are obtained by cylindrical symmetry. 

We may then apply a standard inverse Radon transform (see, for instance, \cite{Nat}) to $R_1g(\tau)$ to recover $g(r,z_0)$ for all $r=|x'|$ and $z_0=x_0^3+t$. Here, $(0,0,z_0)$ corresponds to the intersection point between the beam axis $\mu\mapsto \mu\vec{e}_3$ and the orthogonal (observation) line $\mu\mapsto x_0+t\vec{e}_3 + \mu\theta_0$. 
The inverse Radon transform of rotationally symmetric functions in fact admits the following explicit expression:
\begin{equation}\label{eq:invR}
g(r,z_0) = -\frac{1}{\pi}\int^\infty_r\frac{\frac{d}{d\tau}R_1g(\tau)}{\sqrt{\tau^2-r^2}}d\tau.
\end{equation}

In the presence of several detector arrays modeled by several choices of $x_j\in X$, we obtain from the previous formula the reconstruction of $g(r,x^3)$ for values of $x^3$ sufficiently close to $x^3_j$, $0\leq j\leq J$, and for all $r\geq 0$. See Fig.\ref{fig:geom} for a case with $j=3$ allowing us to reconstruct the beam profile at several positions along the beam axis.






%
\subsection{Determining the beam's main features.}\label{sec:det_beam_feat}
Let us assume that we have access to $g(x',x^3)$ as described above for several values along the profile $z=x^3$. We now propose a reconstruction of the source location (where the beam's source is modeled as a delta function) for different measurement scenarios and prior constraints. 

We assume that the source is $F_0>0$ times a delta function at $x=0$ and $\theta=\vec{e}_3$. From \eqref{meas32}, we thus obtain that $g(x',z)=C_0 e^{-\lambda z} \mathcal{F}^{-1}_{X'}\left[e^{-A_{s}(z)|\xi|^{2s}}\right] \left(x'\right)$ for $z>0$, where $C_0=F_0 \sigma$. Integrating the above expression in $x'\in\Rm^2$ amounts to evaluating the Fourier transform at $\xi=0$ so that $\int_{\Rm^2} g(x',z)dx' = C_0 e^{-\lambda z}$. By computing this quantity for $z_0<z_1$ and then computing its ratio, we reconstruct $e^{-\lambda(z_1-z_0)}$ and hence $\lambda$ since $z_1-z_0$ is known (as the distance between the two detector arrays). Thus, $C_0e^{-\lambda z_0}$ and $\lambda$ are known at this stage, while $C_0$ and $z_0$ remain unknown. 



We now use the reconstructed profile $g(x)$ only at the beam's center $g(z):=g(0,z)$. We assume $g(z)$ known for a number of values of $z\in Z$ as described in the preceding paragraph. The simplest setting is when the set $Z$ is finite. 

From the above considerations, we thus obtain
\begin{equation}\label{eq:g0}
  g(z) = C_0 e^{-\lambda z}\mathcal{F}^{-1}_{X'}\left[e^{-A_{s}(z)|\xi|^{2s}}\right] \left(0\right) =   \frac{C_0 e^{-\lambda z}}{4\pi^2}\int e^{-A_{s}(z)|\xi|^{2s}} d\xi.
\end{equation}

Unlike the ballistic model, the Fermi pencil beam model accounts for beam dispersion which is represented with the last integral above. Moreover, following a computation summarized in the appendix, we observe  that 
\begin{equation}\label{eq:A0}
  A_s(z) = \Big( \dfrac{C_0 e^{-\lambda z}\Gamma(1/s)}{4\pi sg(z)}\Big)^s.
\end{equation}

Since $z$ is discrete as we do not expect to be able to monitor the beam width for all values of $z=x^3$, we assume a constant diffusion coefficients $D$, hence $\tilde D=\frac{1}{2^{2s}}D$ (see \eqref{eq:tilde}). We then verify from the definition in \eqref{meas32} that 
\begin{equation}\label{eq:As_Dctt}
  A_s(z) = \dfrac{\eps^{2s}D}{2s+1} z^{2s+1},
\end{equation}
so that after inversion:
\begin{equation}\label{eq:form1}
z = \left(\frac{2s+1}{\eps^{2s}D}\right)^{\frac{1}{2s+1}}\left(\dfrac{C_0 e^{-\lambda z}\Gamma(1/s)}{4\pi sg(z)}\right)^{\frac{s}{2s+1}}.
\end{equation}
Therefore, if $\eps^{2s}D$ and $s$ are {\em known a priori}, the measurement of $g(z_0)$ for any $z_0$ allows us to reconstruct $z_0$ itself, i.e., the distance from the measured point along the axis to the source location (since $C_0e^{-\lambda z_0}$ is known). Note that $\eps^{2s}D$ is the natural diffusion coefficient appearing in \eqref{u_fFP}.

We now consider a more general setting where $\eps^{2s}D$ and $s$ are also unknown. We may recast the above relations as
\begin{equation}\label{eq:formg}
  g(z) = C_0 e^{-\lambda z}\dfrac{\Gamma(1/s)}{4\pi s} \Big(\dfrac{2s+1}{\eps^{2s}D}\Big)^{\frac1s} z^{-2-\frac1s}.
\end{equation}
Knowledge of $g(z)$ at more than three values of $z$ therefore allows us to reconstruct all coefficients $(\eps^{2s}D,s,z)$ in principle. 

We define $G(t):=g(z_0)/g(z_0+t)$ and find that 
\[
  G(t) = e^{\lambda t}\Big(1+\frac t{z_0}\Big)^{2+\frac1s}.
\]
Knowledge of $G(t_1)$  and $G(t_2)$ for $0<t_1<t_2$ provides a unique reconstruction of the source location $z_0$ and the fraction parameter $s$.  Indeed, define $\alpha=z_0^{-1}$ and $\mu=2+\frac1s$ so that $\ln G(t)-\lambda t =\mu\ln(1+\alpha t)$. We compute
\[
 \partial_\alpha \frac{\ln G(t_1)-\lambda t_1}{\ln G(t_2)-\lambda t_2} = \frac{t_1t_2 [H(t_2)-H(t_1)]}{(1+\alpha t_1)(1+\alpha t_2)\ln^2(1+\alpha t_2)} ,\quad H(t):=(\frac1t+\alpha)\ln(1+\alpha t).
\]
We then obtain $H'(t)=t^{-2}(\alpha t-\ln (1+\alpha t))>0$ so that $\alpha\mapsto \frac{\ln G(t_1)-\lambda t_1}{\ln G(t_2)-\lambda t_2}$ is a strictly increasing function of $\alpha$ when $t_2>t_1>0$. This uniquely determines $\alpha>0$ and hence $\mu$ since $\ln(1+\alpha t_1)$ is now known. Thus $z_0$ and $s$ are uniquely characterized by $g(z_0+t_j)$ for $j=0,1,2$ and $t_0=0<t_1<t_2$. 
It is then straightforward to reconstruct $\eps^{2s}D$ from $g(z)$ in  \eqref{eq:formg} once $s$ and $z=z_0$ are known and $C_0$ from $C_0e^{-\lambda z_0}$ and $\lambda$. 

\medskip

To {\em summarize} the above derivation, we observe that when the beam parameters are constant in $z$, then a finite number of (at least three) measurements of $g(z_j)$ combined with $F_0\sigma e^{-\lambda z_j}=\int_{\Rm^2} g(x',z_j)dx'$ uniquely determine $(\sigma F_0,\eps^{2s}D,\lambda,s,z_0)$. The values of $g(z_j)$ are obtained from the off-axis measurements by an explicit inverse Radon transform \eqref{eq:invR}.  

The parameters $(\eps^{2s}D,\lambda,s)$ characterize the turbulent atmosphere; $z_0$ is a property of the beam, while $C_0=\sigma F_0$ quantifies the strength of the off-axis measurements as a combination of source strength and wide-angle scattering. Only the product $\sigma F_0$ may be reconstructed unambiguously without additional prior information.

\section{Generalizations and remarks} 

\subsection{Errors in line integral measurements}
The inversion procedure presented in the previous section relies on the explicit form of the pencil-beam approximations. The error associated to using such approximate models for the laser beam instead of the more accurate Fokker-Planck may be estimated as follows.

Let $u_S^1$ and $u_S^2$ be the off-axis particle densities solution to \eqref{off_axis_light} with respective source terms $\int \sigma(x,\zeta,\theta)u(x,\zeta)$ and $\int \sigma(x,\zeta,\theta)\fu(x,\zeta)$. We are denoting here by $u$ the Fokker-Planck solution while $\fu$ corresponds to its pencil-beam approximation. Explicitly, we have
\begin{align*}
    u_S^1(x,\theta) & = \int^\infty_0 \int_{\Sm^2}\sigma(x-t\theta,\zeta,\theta)u(x-t\theta,\zeta)d\zeta dt \\
    u_S^2(x,\theta) &= \int^\infty_0 \int_{\Sm^2}\sigma(x-t\theta,\zeta,\theta)\fu(x-t\theta,\zeta)d\zeta dt.
\end{align*}

Given a Lipschitz function $\psi$ with $\|\psi\|_\infty\leq 1$ and $\text{Lip}(\psi)\leq \kappa$ with support contained on a compact set $\omega\subset\Rm^3$ including our off-axis measurements, we consider the solution $\varphi(x,\theta)=\int^\infty_0\psi(x+t\theta,\theta)dt$ of the equation
$
-\theta\cdot\nabla_x\varphi(x,\theta)=\psi
$,
which satisfies $\|\varphi\|_\infty\leq 1$ and $\text{Lip}(\varphi)\leq \kappa$.
Then,
\[
\begin{aligned}
\int \psi(u_S^1-u_S^2)dxd\theta &= \int \varphi(\theta\cdot\nabla_x u_S^1- \theta\cdot\nabla_x u_S^2)dxd\theta \\
&=  \int \sigma(x,\zeta,\theta)\varphi(x,\theta)(u(x,\zeta) - \fu(x,\zeta))d\zeta dxd\theta,
\end{aligned}
\]
which yields
$
\int \psi(u_S^1-u_S^2)dxd\theta \lesssim \mathcal{W}^1_\kappa(u,\fu)
$
for all $\psi$ as above. Taking supremum among all those $\psi$ and recalling the results in Theorem \ref{thm:W_est}, we obtain the estimate
\begin{equation}\label{eq:errormeas}
\mathcal{W}^1_{\kappa,\omega}(u_S^1,u_S^2)\leq C\kappa^{s'} \epsilon^{2s'}
\end{equation}
for $s'$ and the constant $C$ as in the theorem (where $s'=1$ if $s=1$) and where $\mathcal{W}^1_{\kappa,\omega}$ is the $(1,\kappa)$-Wasserstein distance restricted to $\omega$.

We thus obtain that the measurement errors generated by replacing Fokker-Planck models by their Fermi pencil beam approximation are small in the above sense when $\eps$ is small. 

We also refer to \cite{BJ} for the effect of such measurement errors on the reconstruction of the function $g(z)$ from its line integrals.

\subsection{The local case $s=1$} 

Recalling definition \ref{def:Fpb} in the local case $s=1$, the Fermi pencil-beam takes the form of the following Gaussian beam:
\[
\begin{aligned}
U(X,V) &= e^{-\int^{X^3}_0\tilde{\lambda}(r)dr}\mathcal{F}^{-1}_{X',V}\left[ e^{-\int^{X^3}_0|\eta+(X^3-t)\xi|^{2}D(t)dt}\right]\\
&=\frac{e^{-\int^{X^3}_0\tilde{\lambda}(r)dr}}{(4\pi)^2(E_2E_0-E_1^2)}\exp\left\{-\frac{E_0|X'-X^3V|^2+2E_1(X'-X^3V)\cdot V+E_2|V|^2}{4(E_2E_0-E_1^2)}\right\}
\end{aligned}
\]
with $E_k(X^3):= \tilde D\int^{X^3}_0t^k dt$ and $2^{2}\tilde D=D$ here. 
We use this exponentially decaying function to define the pencil-beam approximation in \eqref{Fpb}.

The reconstruction methodology employed in section \ref{sec:rec} applies to any $s\in(0,1]$. 
However, the determination of  $(\sigma F_0,\epsilon^2D,z,\lambda)$ simplifies when $s=1$. In particular, \eqref{eq:g0} becomes
\begin{equation}\label{g_s1_rad}
g(z)=\frac{C_0e^{-\lambda z}}{4\pi A_{1}(z)},
\end{equation}
in a neighborhood of some $z_0>0$. The factor $C_0e^{-\lambda z}$ as well as $\lambda$ are obtained from the integral $\int_{\Rm^2}g(x',z)dx'$ as in  section \ref{sec:det_beam_feat}. This yields $A_1(z)=\frac{\epsilon^{2} D}{3}z^{3}$.
For a given $t>0$, we compute $A_1(z+t)/A_1(z)=(1+t/z)^3$ and then obtain the distance $z$ following the relation
\[
z = \frac{t}{\left(\frac{A_1(z+t)}{A_1(z)}\right)^{1/3}-1}.
\]
Finally, from knowledge of $z$ and  $A_1(z)$ we easily determine the factor $\epsilon^2D$.

\medskip


\subsection{Broader source terms}

Several generalizations of the above reconstructions may be considered in the presence of additional measurements. For instance, if measurements are available on a continuum of values of $z_0$ (with $\Sigma$ involving a continuum of values of $x_0$), then the parameters $(\lambda(z), D(z))$ and possibly $s$ as well may be allowed to vary in the $z-$variable; we do not consider this particular setting any further.

Here, we consider a generalization with $D$ and $\lambda$ still constant but with a laser source that is spatially broad but still narrow in its direction of emission.
The laser beam is modeled as a solution to \eqref{fFPB} with an unknown source
%
\[
U(X',0,V) = G(X',V):=h(X')\delta(V),
\]
for a nonnegative and integrable function $h$, supported inside the ball of radius $\rho$ with $\rho+\epsilon\ll L_0$ (the latter parameter represents the size of the camera's screen).

In the Fourier domain this solution takes the form
\[
\mathcal{F}_{X',V}[U](\xi,X^3,\eta) = e^{-\lambda X^3-\frac{D}{2^{2s}}\int^{X^3}_0|\eta+(X^n-t)\xi|^{2s}dt}\mathcal{F}_{X',V}[G](\xi,\eta+X^3\xi).
\]
The available measurements thus satisfy
\[
m(x,\theta) = \int^\infty_0g(x+t\theta)dt +  O(\epsilon^{2s'}),
\]
for $s'\in(0,s)$ close to $s$ (with the error corresponding to the relative one), and with $g(x)$ given by
\[
g(x)=\sigma e^{-\lambda x^3}\mathcal{F}^{-1}_{X'}\left[e^{-A_s(x^3)|\xi|^{2s}}\mathcal{F}_{X',V}[G](2\epsilon\xi,2\epsilon \xi x^3)\right](x').
\]
Taking into account the explicit form of $G$, the above reduces to
\[
g(x)=\sigma e^{-\lambda x^3}\mathcal{F}^{-1}_{X'}\left[e^{-A_s(x^3)\xi^{2s}}\mathcal{F}_{X'}[h](2\epsilon\xi)\right](x').
\]
In the simpler case of a cylindrical beam, that is, with $h(X')=h(|X'|)$, the radial symmetry of $h$ is inherited by $g$ and therefore, for each $t\ll 1$, a single measurement of the form $m(x_0+t\phi_0+\tau \psi_0,\theta_0)$ with $\phi_0,\varphi_0$ and $\theta_0$ as in \ref{sec:beam_meas} and $t^2+\tau^2<L_0^2$, gives us all the line integrals on the plane passing through $x_0+t\phi_0$ and perpendicular to $\phi_0=\vec{e}_3$.

In the more general case of $h$ not necessarily radial-symmetric, we are forced to acquire a larger set of measurements. The tomographic procedure introduced above involving the inversion of a Radon transform applies in this case provided we observe the beam from an array of multiple cameras totally or partially surrounding the beam. 

\begin{figure}
\centerline{
\includegraphics[scale=0.4]{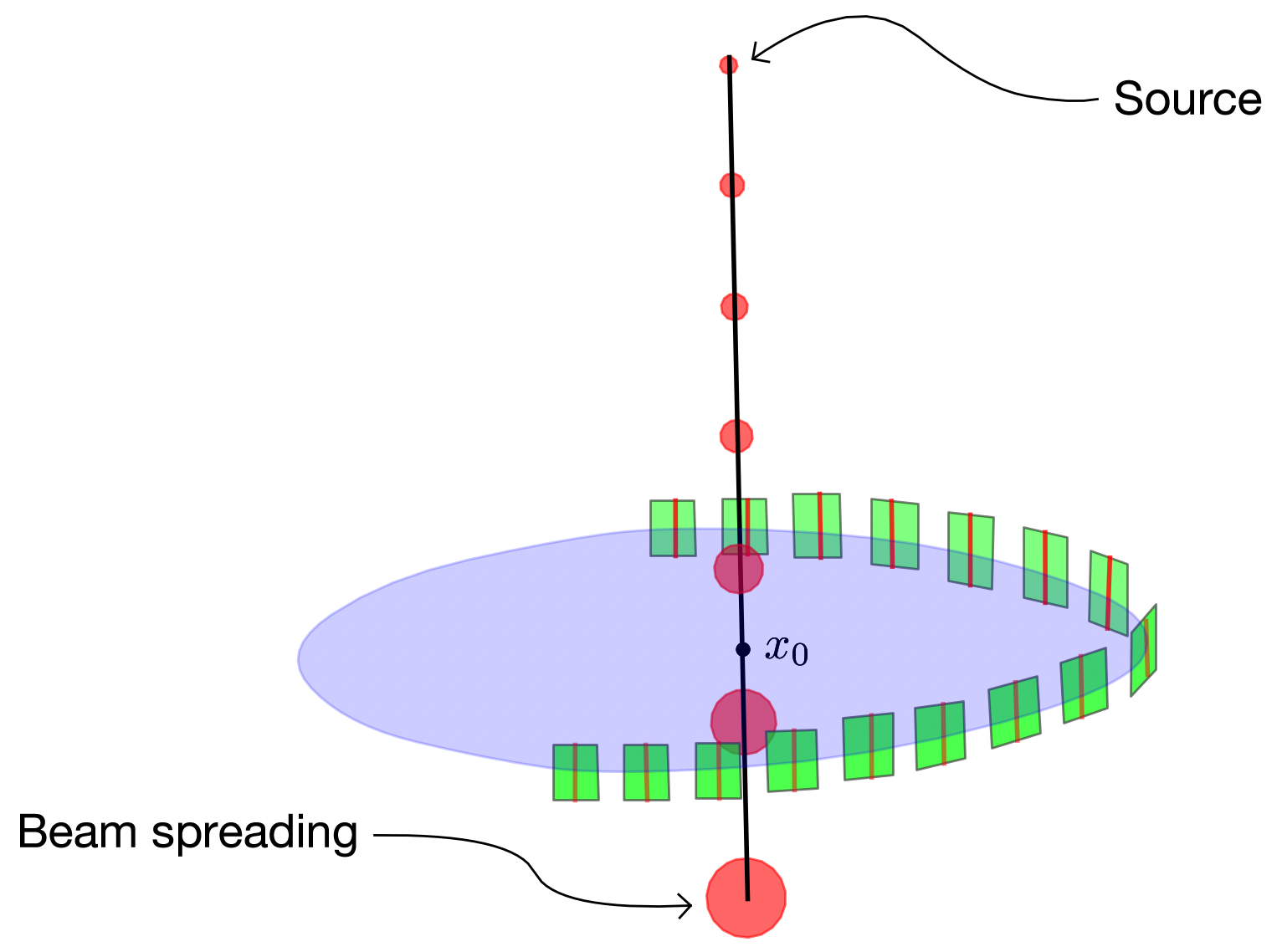}
}
\caption{Tomographic measuring geometry for non-cylindrical beams. In the presence of cylindrical symmetry one detector is enough.}
\label{fig:noncylidrical}
\end{figure}

An array of multiple cameras placed on a plane perpendicular to the beam's axis can in principle measure all the line integrals passing near the intersection point of the axis and the plane (say $x_0=(0,0,z_0)$). Therefore, a Radon transform inversion leads to the determination of $g(x',z_0)$ for all $x'\in\Rm^2$. See figure \ref{fig:noncylidrical} for an schematic of this measurement geometry.
In addition, by Fourier transforming $g(x)$ with respect to $x'$ we have access to
\begin{equation}\label{Fg}
\mathcal{F}_{x'}[g](\xi,z_0)=\sigma e^{-\lambda z_0}e^{-A_s(z_0)|\xi|^{2s}}
\mathcal{F}_{X',V}[h](2\epsilon\xi),\quad\forall \xi\in\Rm^2,
\end{equation}
from which we realize that 
\[
\begin{aligned}
\mathcal{F}_{x'}[g](0,z_0)=\int_{\Rm^2} g(x',z_0)dx'
&=\sigma e^{-\lambda z_0}\int_{\Rm^2}h(x')dx'.
\end{aligned}
\]
Consequently, and denoting $C_0=\sigma \int_{\Rm^2}h(x')dx'$, we obtain the quantity $C_0e^{-\lambda z_0}$. Computing this for both, $z_0<z_1,$ we obtain $\lambda$ (as done previously).

For a fixed (and known) $t>0$, let us now define $\mathcal{G}(\xi,t)= \mathcal{F}_{x'}[g](\xi,z_0)/\mathcal{F}_{x'}[g](\xi,z_0+t)$ and observe that 
\[
\ln\left(\mathcal{G}(\xi,t)\right) - \lambda t =  |\xi|^{2s}\left(A_s(z_0+t)-A_s(z_0)\right).
\]
Note that while $\mathcal{F}_{x'}[g](\xi,z_0)$ depends on $\xi$, the ratio $\mathcal{G}(\xi,t)$ depends only on $|\xi|$.
Evaluating the above at $|\xi|=1$ provides the difference $A_s(z_0+t)-A_s(z_0)$. Subsequently choosing $|\xi|=e$, we apply logarithm to the previous expression and deduce that
\[
s = \frac{1}{2}\ln\left(\frac{\ln(\mathcal{G}(\xi,t)-\lambda t)}{A_s(z_0+t)-A_s(z_0)}\right).
\]
In homogeneous media (i.e., with constant $D$ and $\lambda$), we use the explicit form of $A_s(z_0)$ in \eqref{eq:As_Dctt} to get
\begin{equation}\label{logG}
\ln\left(\mathcal{G}(\xi,t)\right) - \lambda t = \frac{|\xi|^{2s}\epsilon^{2s}D}{(2s+1)(z_0)^{2s+1}}\left(\left(1 +\frac{t}{z_0}\right)^{2s+1}-1\right).
\end{equation}
Then, using the above with $t_1<t_2$ and denoting $\alpha=z_0^{-1}$, we have
\[
H(\alpha)=\frac{\ln\left(\mathcal{G}(\xi,t_1)\right) - \lambda t_1}{\ln\left(\mathcal{G}(\xi,t_2)\right) - \lambda t_2} = \frac{\left(1 +\alpha t_1\right)^{2s+1}-1}{\left(1 +\alpha t_2 \right)^{2s+1}-1}.
\]
This is a strictly increasing function of $\alpha$ as we can verify by computing its derivative. Indeed,
\[
\begin{aligned}
\partial_\alpha H&=(2s+1)\frac{t_1(1+\alpha t_1)^{2s}(\left(1 +\alpha t_2 \right)^{2s+1}-1)-t_2(1+\alpha t_2)^{2s}(\left(1 +\alpha t_1 \right)^{2s+1}-1)}{(\left(1 +\alpha t_2 \right)^{2s+1}-1)^2}\\
&=\frac{t_1t_2}{(\left(1 +\alpha t_2 \right)^{2s+1}-1)^2}\int^\alpha_0\left((1+\alpha t_1)^{2s}(1 +r t_2 )^{2s}-(1+\alpha t_2)^{2s}(1 +r t_1 )^{2s}\right)dr,\\
\end{aligned}
\]
with a negative integrand since the function $r\mapsto\frac{1+r t_1}{1+r t_2}$ is strictly decreasing.

We then conclude that $\alpha = z_0^{-1}$ is uniquely determined by the quantity $H(\alpha)$ and hence so is $s$. Subsequently, we can determine $C_0=\sigma\left(\int_{\Rm^2}h(x')dx'\right)$, and $\epsilon^{2s}D$ (and then $A_s(z_0)$) for instance from \eqref{logG}. 

Lastly, going back to \eqref{Fg}, we are able to determine a rescaled version of $h$ up to a constant factor given by $\left(\int_{\Rm^2}h(x')dx'\right)^{-1}$. This follows from the relation
\[
\left(\int_{\Rm^2}h(x')dx'\right)^{-1}\mathcal{F}_{x'}[h](2\epsilon \xi) = \frac{\mathcal{F}_{x'}[g](\xi,z_0)}{C_0 e^{-\lambda z_0}e^{-A_s(z_0)|\xi|^{2s}}},
\]
so that Fourier transforming both sides yields
\[
\left(\int_{\Rm^2}h(x')dx'\right)^{-1}\frac{1}{(2\epsilon)^2}h\left(\frac{x'}{2\epsilon}\right) = \frac{\mathcal{F}_{x'}^{-1}\left[\mathcal{F}_{x'}[g](\xi,z_0)e^{A_s(z_0)|\xi|^{2s}}\right](x')}{C_0 e^{-\lambda z_0}}.
\]
Recalling the relation between the stretched and macroscopic variables: $X'=x'/2\epsilon$, we determine from the above (up to a constant factor), the source function $h$ in its natural (stretched) coordinates, namely,
\[
\frac{h_\epsilon(X')}{\int_{\Rm^2} h_\epsilon(Y)dY} = \frac{\mathcal{F}_{x'}^{-1}\left[\mathcal{F}_{x'}[g](\xi,z_0)e^{A_s(z_0)|\xi|^{2s}}\right](2\epsilon X')}{C_0 e^{-\lambda z_0}}.
\]

\section{Conclusions}\label{sec:conclu}
The reconstruction of the profile of a laser beam propagating through a turbulent atmosphere from limited off-axis measurements is a difficult task.  We propose here a linear macroscopic description of the beam spreading based on a (possibly fractional) Fermi pencil beam equation. Such approximations are accurate in the regime of small mean-free path large transport-mean-free path consistent with the narrow laser beam hypothesis. Moreover, by neglecting back-scattering, they admit sufficiently explicit expressions that are amenable to parameter inversions. 

The off-axis measurement assumptions are as follows. We assume the presence of detector arrays away from the path of the beam. Light detection is modeled as wide-angle single scattering off of the laser beam. The scattering amplitude and the laser source amplitude remain unknown but it is assumed that their combined effect is large enough that it can be detected. 

At least two detector arrays are necessary to triangulate the line segment along which the beam propagates. Then, one sufficiently large detector, or several detectors along the path, allow us to evaluate the whole beam structure under suitable assumptions. The explicit influence of the parameters of the model on beam spreading enables us to reconstruct them. Such macroscopic parameters include the source location and the main features of the turbulence through which it propagates (diffusion coefficient and fractional power). Note that a model based on expansions of radiative transfer solutions into successive scattering events as in \cite{C-SPIE-03} would model the beam intensity as a ballistic component, which cannot possibly provide information on the source location, say, since beam spreading is absent from the model.

Our explicit reconstructions are based on the inversion of the laser intensity at the center of a radially symmetric beam at several locations along its path. The latter intensity may be estimated by applying a standard inverse Radon transform on the available detector array measurements. Other, possibly more stable, reconstruction procedures from the same available data are certainly possible. Our results provide proof of concept that (fractional) Fermi pencil beam models allow for the reconstruction of macroscopic laser beam features including what was our main motivation here: their source location.

\section*{Acknowledgment}
This research was partially supported by the Office of Naval Research, Grant N00014-17-1-2096 and by the National Science Foundation, Grant DMS-1908736.

Most of the work presented in this article was done while BP was a W. H. Kruskal Instructor at the University of Chicago. BP would like to thank the University of Chicago and in particular the Department of Statistics for their hospitality and generosity throughout those years. 

\begin{appendix}
\section{}
To compute $ \int_{\Rm^3} e^{-A_{s}|\xi|^{2s}} d\xi$ we consider polar coordinates and obtain
\[
\int_{\Rm^3} e^{-A_{s}|\xi|^{2s}} d\xi= 2\pi\int_0^\infty \rho e^{-A_{s}\rho^{2s}} d\rho.
\]
Recalling the definition of the Gamma function, $\Gamma(z)=\int^\infty_0 x^{z-1}e^{-x}dx$, defined for $\text{Re}(z)>0$, we are able to recast the integral as
\[
\int_{\Rm^3} e^{-A_{s}|\xi|^{2s}} d\xi= 2\pi A_{s}^{-1/s}\int^\infty_0 te^{-t^{2s}}dt = \frac{\pi }{sA_{s}^{1/s} }\int^\infty_0 t^{1/s-1}e^{-t}dt = \frac{\pi }{sA_{s}^{1/s}}\Gamma\left(\frac{1}{s}\right).
\]
\end{appendix}


\end{document}